\documentclass[11pt, reqno]{amsart}
\usepackage[margin=3.6cm]{geometry}

\usepackage{amsmath, amsfonts, amsthm, amssymb, dsfont}
\usepackage{exscale}
\usepackage{cite}
\usepackage{epsfig}
\usepackage{amscd}
\usepackage{graphics}
\usepackage{color}

\renewcommand{\geq}{\geqslant}
\renewcommand{\leq}{\leqslant}

\newtheorem{theorem}{Theorem}

\newtheorem{proposition}{Proposition}[section]

\newtheorem*{main-theorem}{Main Theorem}
\newtheorem*{theorem*}{Theorem}

\theoremstyle{definition}

\newtheorem{remark}[proposition]{Remark}

\newtheorem*{remark*}{Remark}

\numberwithin{equation}{section}

\def\phi{\varphi}

\def\reals{{\mathbb R}}

\def\O{{\mathcal O}}

\def\phi{\varphi}

\def\be{\begin{eqnarray*}}
\def\ee{\end{eqnarray*}}
\def\ben{\begin{eqnarray}}
\def\een{\end{eqnarray}}

\def\L2R{L_{\text{Rest}}^2}

\def\11{\mathds{1}}

\def\L2c{L^2_{\text{comp}}}

\def\p{\partial}

\def\bu{\bar{u}}

\begin{document}

\title[Neumann data on triangles]{Equidistribution of 
  Neumann data mass on triangles}

%% \author[Y. Canzani]{Yaiza Canzani}
%% \address[Y. Canzani]{Department of Mathematics, Harvard University.\medskip}
%%  \email{canzani@math.harvard.edu}

\author[H. Christianson]{Hans Christianson}
\address[H. Christianson]{ Department of Mathematics, University of North Carolina.\medskip}
 \email{hans@math.unc.edu}

%% \author[J. Toth]{John  A. Toth}
%% \address[J. Toth]{Department of Mathematics and Statistics, McGill
%%   University.\medskip} 
%% \email{jtoth@math.mcgill.ca}

%% \author[J. Galkowski]{Jeffrey Galkowski}
%% \address[J. Galkowski]{Department of Mathematics and Statistics,
%%   McGill University.\medskip}
%% \email{}

% \thanks{}
%
%\subjclass[2000]{}
%\keywords{}
%
%
\begin{abstract}
  In this paper we study the behaviour of the Neumann data 
of Dirichlet eigenfunctions on triangles.   We prove that the  $L^2$
norm of the (semi-classical) Neumann data on each side is equal to the
length of the side divided by the area of the triangle.  The novel
feature of this result is that it is {\it not} an asymptotic, but an
exact formula.  The proof is by simple integrations by parts.

\end{abstract}

\maketitle

\section{Introduction}
Given a compact surface or manifold with boundary,   it is an
interesting question to consider restrictions of 
eigenfunctions to hypersurfaces; either the Dirichlet data or Neumann
data (or both, the Cauchy data) can be considered.  Perhaps the
simplest question is to consider boundary values.  That is, if we
consider Dirichlet (respectively Neumann) eigenfunctions, we may try
to study the Neumann (respectively Dirichlet) data on the boundary.

In this short note, we consider one of the simplest planar domains, a planar triangle $T$.  Our main
result is that the $L^2$ mass of the semi-classical Neumann data on each side of $T$ {\it equals} the
length of the side divided by the area of $T$.  It should be
emphasized that  these formulae are equalities, not  asymptotics or
estimates.

\begin{theorem}
  \label{T:main}
Let $T$ be a planar triangle  with sides $A,B,C$, of lengths $a,b,c$
respectively.  Consider the (semi-classical) Dirichlet eigenfunction
problem:
\begin{equation}
  \label{E:efcn}
\begin{cases}
  (-h^2 \Delta -1)u = 0 , \text{ in } T, \\
  u |_{\p T} = 0,
\end{cases}
\end{equation}
and assume the eigenfunctions are normalized
$
  \| u \|_{L^2(T)} = 1.
$

  Then the (semi-classical) Neumann data on the boundary satisfies
  \[
  \int_A | h \p_\nu u |^2 dS = \frac{a}{\text{Area}(T)},
  \]  
 \[
  \int_B | h \p_\nu u |^2 dS = \frac{b}{\text{Area}(T)},
  \]
and
 \[
  \int_C | h \p_\nu u |^2 dS = \frac{c}{\text{Area}(T)},
  \]
where $h \p_\nu$ is the semi-classical normal derivative on $\p T$,
$dS$ is the arclength measure, and $\text{Area}(T)$ is the area of the
triangle $T$.

\end{theorem}

We pause to note briefly that the 
 the semiclassical parameter $h$ takes discrete values as $h \to
0$ (reciprocals of eigenvalues).

\begin{remark}

  We are calling this ``equidistribution'' of Neumann mass since it
  says that the Neumann data has the same mass to length ratio.  Of
  course it does not say anything about local equidistribution to
  subsets of the sides.

 To the author's knowledge, no exact formula such as  this exists in
 the previous literature, except in cases where explicit formulae for
 the eigenfunctions are known.  Even in these cases, the formulae
 typically depend on $h$.   A statement such as Theorem \ref{T:main} is false in general for other planar polygons.  See
 Section \ref{S:example} for the example of a square.

   In order to better understand these formulae, in subsequent works, the author will study non-Euclidean triangles,
and higher dimensional problems, as
well as explore weaker lower bounds and interior lower bounds in some
simple polygons.
\end{remark}

\subsection{History}

  Previous results on restrictions primarily focused on upper bounds.
  In the paper of Burq-G\'erard-Tzvetkov \cite{BGT-erest}, restrictions of
  the Dirichlet data to arbitrary hypersurfaces were considered.  An upper
  bound of the norm (squared) of the restrictions of $\O(h^{-1/2})$ was
  proved, and shown to be sharp.  Of course this shows that there are
  {\it some} eigenfunctions with a known lower bound on the norms of restrictions.  In
  the author's paper with Hassell-Toth \cite{CHT-ND}, an upper bound of $\O(1)$ was
  proved for (semi-classical) Neumann data restricted to arbitrary hypersurfaces,  and also shown to be sharp.  Again, this gives a lower and
  upper bound for {\it some} eigenfunctions.

  In the case of quantum ergodic eigenfunctions, a little more is
  known.  In the papers of G\'erard-Leichtnam \cite{GeLe-qe} and Hassell-Zelditch \cite{HaZe}, the Neumann (respectively
  Dirichlet) boundary data of Dirichlet (respectively Neumann) quantum ergodic eigenfunctions is
  studied, and shown to have an asymptotic formula for a {\it density
    one} subsequence.  That means that there is a lower bound, and
  explicit local asymptotic formula in this special case, at least for
  most of the eigenfunctions.  Similar statements were proved for
  interior hypersurfaces by Toth-Zelditch \cite{ToZe-1, ToZe-2}.
  Again, potentially a sparse subsequence may behave differently.  In
  the the author's paper with Toth-Zelditch \cite{CTZ-1}, an asymptotic formula for the
  whole weighted Cauchy data is proved for the entire sequence of
  quantum ergodic eigenfunctions, however it is impossible to separate
  the behaviour of the Dirichlet versus Neumann data.

\section{Proof of Theorem \ref{T:main}}

Assume the sides $A,B,C$ are listed in clockwise
orientation.
We assume that $A$ is the shortest side, followed by $B$
and $C$ with respective lengths  $a\leq b \leq c$.

We use rectangular coordinates $(x,y)$ in the plane, and  orient our
triangle so that the corner between $B$ and $C$ is at the origin
$(0,0)$.  
  We further assume that the side $A$ is parallel to
the $y$ axis.

\begin{figure}
\hfill
\centerline{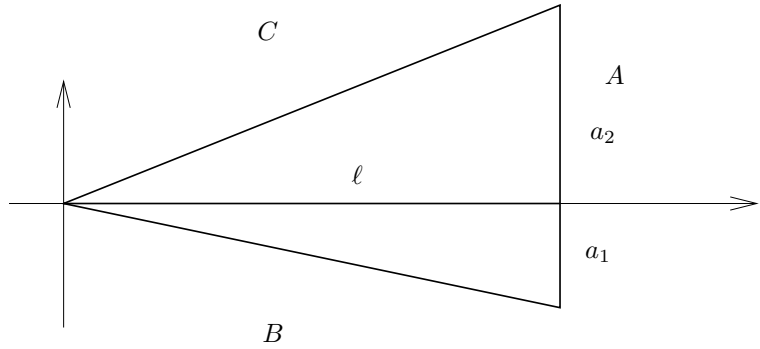}
\caption{\label{F:acute}  Setup for acute (and right) triangles}
\hfill
\end{figure}

\begin{figure}
\hfill
\centerline{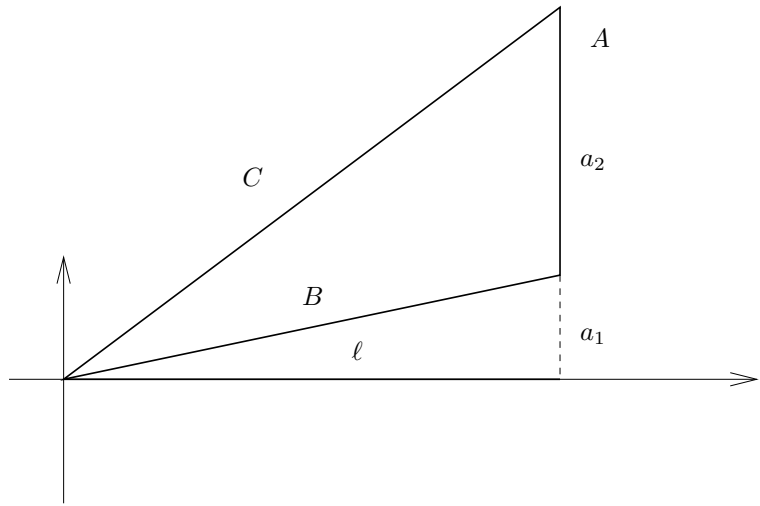}
\caption{\label{F:obtuse}  Setup for obtuse triangles}
\hfill
\end{figure}

We break our analysis into the two cases of acute triangles (including
right triangles) and obtuse.  See Figures \ref{F:acute} and
\ref{F:obtuse} for a picture of the setup.

\subsection{Acute triangles}
Let $\ell$ be the segment on the $x$ axis beginning at $(0,0)$ and
perpendicular to the side $A$.  Of course $\ell$ can be computed in
terms of the sides, but its value is not necessary for this
computation, other than to note that the area of $T$ is $a \ell /2$.
Write $A = A_1 \cup A_2$, where $A_1$ is the part of $A$ under the $x$
axis and $A_2$ is the part above.  Let $a_1, a_2$ denote the
respective sidelengths.

We can parametrize $B$ and $C$ with respect to $x$.
\[
C = \left\{ (x,y) \in \reals^2 : y = \frac{a_2}{\ell} x, \,\, 0 \leq x \leq
\ell \right\},
\]
and
\[
B = \left\{ (x,y) \in \reals^2 : y = - \frac{a_1}{\ell} x, \,\, 0 \leq x
\leq \ell \right\}.
\]

Then the arclength parameters are
\[
\gamma_C = ( 1 + (a_2 / \ell)^2 )^{1/2} = \frac{ (\ell^2 + a_2^2
  )^{1/2} }{\ell} = \frac{c}{\ell},
\]
and
\[
\gamma_B = ( 1 + (a_1 / \ell)^2 )^{1/2} = \frac{ (\ell^2 + a_1^2
  )^{1/2} }{\ell} = \frac{b}{\ell},
\]
 and the unit tangent vectors are
\[
\tau_C = \left(1 , \frac{ a_2}{\ell} \right) \gamma_C^{-1} = \left(\frac{\ell}{c} ,
\frac{a_2}{c}\right)
\]
and
\[
\tau_B = \left(1 , -\frac{ a_1}{\ell} \right) \gamma_B^{-1} = \left(\frac{\ell}{b} ,
-\frac{a_1}{b}\right).
\]

From this we have the outward unit normal vectors
\[
\nu_C = \left(- 
\frac{a_2}{c}, \frac{\ell}{c} \right)
\]
and
\[
\nu_B = \left(- 
\frac{a_1}{b}, -\frac{\ell}{b} \right).
\]
Of course the outward normal to $A$ is $\nu_A = (1,0)$.

We are assuming Dirichlet boundary conditions, which implies that the
tangential derivatives of $u$ vanish on $\p T$.  That is,
\[
\p_y u = 0
\]
on $A$, and
\[
\tau_C \cdot \nabla u = \frac{\ell}{c} \p_xu   + 
\frac{a_2}{c} \p_y u = 0
\]
on $C$.
Similarly,
\[
\tau_B \cdot \nabla u = \frac{\ell}{b} \p_xu   - 
\frac{a_1}{b} \p_y u = 0
\]
along $B$.  Rearranging, we have
\[
h\p_x u = - \frac{a_2}{\ell} h \p_y u
\]
on $C$ and
\[
h\p_x u =  \frac{a_1}{\ell} h \p_y u
\]
on $B$.

Making the substitutions, along $C$ we have
\begin{align*}
  h\p_{\nu_C} u & = \nu_C \cdot h\nabla u \\
  & = - 
  \frac{a_2}{c} h \p_x u + \frac{\ell}{c} h \p_y u \\
  & = \left(\frac{a_2^2}{c \ell} h \p_x u + \frac{\ell}{c}\right) h \p_y u \\
  & = \left( \frac{a_2^2 + \ell^2}{c \ell} \right) h \p_y u \\
  & = \frac{c}{\ell} h \p_y u.
\end{align*}
Hence
\[
h \p_y u = \frac{\ell}{c} h \p_{\nu_C} u
\]
on $C$.  Substituting again, we have
\[
h \p_x u = - \frac{a_2}{\ell} h \p_y u = - \frac{a_2}{c} h \p_{\nu_C}
u
\]
along $C$.
Similarly, along $B$, we have
\[
h \p_y u = - \frac{\ell}{b} h \p_{\nu_B} u
\]
and
\[
h \p_x u = - \frac{a_1}{b} h \p_{\nu_B}
u.
\]

We now consider the vector field
\[
X = (x + m) \p_x + (y + n ) \p_y,
\]
where $m,n$ are parameters independent of $x$ and $y$.  Since $m \p_x$ commutes with $-h^2 \Delta$
as well as $n \p_y$, the usual computation yields
\[
  [-h^2 \Delta -1 , X ] = -2h^2 \Delta.
  \]
  Then using eigenfunction equation \eqref{E:efcn}, we have
  \begin{align*}
  \int_T ([-h^2 \Delta -1 , X ] u) \bu dV & = -2\int_T (h^2 \Delta u)
  \bu dV \\
  & = \int_T 2 | u |^2 dV \\
  & = 2,
  \end{align*}
  since $u$ is normalized.

  On the other hand, again using the eigenfunction equation
  \eqref{E:efcn} again, we have
  \begin{align*}
    \int_T & ([-h^2 \Delta -1 , X ] u) \bu dV \\
    & = 
  \int_T ((-h^2 \Delta -1 ) X  u) \bu dV - 
  \int_T (X (-h^2 \Delta -1 ) u) \bu dV \\
  & = \int_T ((-h^2 \Delta -1 ) X  u) \bu dV.
  \end{align*}
  Integrating by parts and using the eigenfunction equation and Dirichlet boundary conditions, we have
  \begin{align*}
    \int_T & ((-h^2 \Delta -1 ) X  u) \bu dV \\
    & = \int_T  ( X  u) ((-h^2 \Delta -1 ))\bu dV \\
& \quad - \int_{\p T}  ( h \p_\nu hX  u) \bu dS  +
    \int_{\p T}  ( hX  u) ( h \p_\nu \bu) dS\\
    & = \int_{\p T}  ( hX  u) ( h \p_\nu \bu) dS,
  \end{align*}
  Hence we have computed:
  \[
  2 = \int_{\p T}  ( hX  u) ( h \p_\nu \bu) dS.
  \]

Let us break up the analysis into the three different sides.  In order
to simplify notation somewhat, set
\[
I_A = \int_A | h \p_\nu u |^2 dS,
\]
and similarly for $B$ and $C$.  Notice we have left the surface
measure $dS$ alone, even though we could write it explicitly in terms
of the arclength parameters computed above.  This is not necessary for
the analysis.

  Returning now to our computations of the normal derivatives, we have
  \begin{align*}
    \int_A &(hX u ) ( h \p_\nu \bu ) dS \\
    & = 
  \int_A (((x + m ) h \p_x + (y + n ) h \p_y )u ) ( h \p_\nu \bu ) dS
  \\
  & = (\ell + m ) I_A,
  \end{align*}
  since $x = \ell$ on $A$.

  Continuing, using that along $C$, we have $y = (a_2/\ell) x$:
  \begin{align*}
    \int_C & (hX u ) (h \p_\nu \bu) dS \\
    & =  \int_C (((x + m ) h \p_x + (y + n ) h \p_y )u ) ( h \p_\nu
    \bu ) dS \\
    & = 
 \int_C \left(((x + m ) h \p_x + \left(\frac{a_2}{\ell} x + n \right) h \p_y )u \right) ( h \p_\nu \bu ) dS
\\& = \int_C \left( \left( (x + m ) \left(- \frac{a_2}{c}\right) + \left( \frac{a_2}{\ell}x + n
\right)\left(\frac{\ell}{c} \right) \right) h \p_{\nu_C} u \right) ( h \p_{\nu_C} \bu ) dS \\
& = \int_C \left(   \left(- \frac{a_2}{c}m +  
\frac{\ell}{c} n\right)  h \p_{\nu_C} u \right) ( h \p_{\nu_C} \bu ) dS \\
& = 
\left(- \frac{a_2}{c}m +  
\frac{\ell}{c} n \right) I_C.
\end{align*}

  Similarly, along $B$, we have $y = -(a_1/\ell) x$:
  \begin{align*}
    \int_B & (hX u ) (h \p_\nu \bu) dS \\
    & =  \int_B (((x + m ) h \p_x + (y + n ) h \p_y )u ) ( h \p_\nu
    \bu ) dS \\
    & = 
 \int_B \left(((x + m ) h \p_x + \left(-\frac{a_1}{\ell} x + n \right) h \p_y )u \right) ( h \p_\nu \bu ) dS
 \\
 & = \int_B \left(\left( (x + m ) \left(- \frac{a_1}{b}\right) + \left( -\frac{a_1}{\ell}x + n
\right)\left(-\frac{\ell}{b} \right) \right) h \p_{\nu_B} u \right) ( h \p_{\nu_B} \bu ) dS \\
& = \int_B \left(   \left(- \frac{a_1}{b}m -  
\frac{\ell}{b} n \right)  h \p_{\nu_B} u \right) ( h \p_{\nu_B} \bu ) dS \\
& = 
\left(- \frac{a_1}{b}m -  
\frac{\ell}{b} n \right) I_B.
\end{align*}

  Summing up, we have:
  \begin{equation}
    \label{E:master-eqn}
  2 = (\ell + m ) I_A + \left(- \frac{a_1}{b}m -  
\frac{\ell}{b} n \right) I_B
+ \left(- \frac{a_2}{c}m +  
\frac{\ell}{c} n \right) I_C.
\end{equation}

First, set $m = n = 0$.  Then we get
$
2 = \ell I_A,
$
or 
\begin{align*}
  I_A & = \frac{2}{\ell} \\
  & = \frac{a}{a \ell /2}  \\
  & = \frac{a}{\text{Area}(T)}.
\end{align*}

Now we observe that the left hand side of \eqref{E:master-eqn} is independent of $m$ and
$n$ so we can differentiate with respect to $m$ and $n$ to get two new
equations.  That is, differentiating both sides of
\eqref{E:master-eqn} with respect to $m$ yields
\[
0 = I_A - \frac{a_1}{b} I_B - \frac{a_2}{c}I_C,
\]
and plugging in the value of $I_A$,  we have
\begin{equation}
  \label{E:d-m}
\frac{a_1}{b} I_B + \frac{a_2}{c}I_C = \frac{2}{\ell}.
\end{equation}
Now differentiating \eqref{E:master-eqn} with respect to $n$, we have
\[
  0 =   -  
\frac{\ell}{b}  I_B
 +  
 \frac{\ell}{c}  I_C,
 \]
 so that
 \[
 I_B = \frac{b}{c} I_C.
 \]
 Plugging in to \eqref{E:d-m}, we have
 \begin{align*}
 \frac{2}{\ell} & = \left( \left(\frac{a_1}{b}\right)\left(\frac{b}{c}\right) + \frac{a_2}{c}\right)
 I_C \\
 & = \left( \frac{a_1}{c} + \frac{a_2}{c} \right) I_C \\
 & = \frac{a}{c} I_C.
 \end{align*}
 Hence
 \[
 I_C = \frac{2c}{a \ell} = \frac{c}{\text{Area}(T)},
 \]
 and back substituting,
 \[
 I_B = \frac{b}{c} I_C = \frac{b}{\text{Area}(T)}.
 \]
 This proves the theorem for acute and right triangles.

 \subsection{Obtuse triangles}
 The proof is nearly the same, with several sign changes.  Using the
 setup in Figure \ref{F:obtuse}, we have
 \[
 C = \left\{ (x,y) : 0 \leq x \leq \ell, \text{ and } y = \frac{(a +
   a_1)}{\ell x} \right\},
 \]
 and
 \[
 B = \left\{ (x, y) : 0 \leq x \leq \ell, \text{ and } y = \frac{a_1}{\ell}
 x \right\}.
 \]
 Similar computations as above lead to the following: along $C$,
 \[
 h \p_y u = \frac{\ell}{c} h \p_{\nu_C} u,
 \]
 and
 \[
 h \p_x u = - \frac{a + a_1}{c} h \p_{\nu_C}.
 \]
 Along $B$ we have
 \[
 h \p_y u = - \frac{ \ell}{b} h \p_{\nu_B}u
 \]
 and
 \[
 h \p_x u =   \frac{a_1}{b} h \p_{\nu_B} u.
 \]

 The same commutator computation holds, and similar substitutions as
 in the acute case yield the equation
 \begin{equation}
   \label{E:master-2}
   2 = ( \ell + m ) I_A + \left(\frac{a_1}{b}m - \frac{\ell}{b}n \right) I_B +
   \left(-\frac{a + a_1}{c} m + \frac{ \ell}{c} n \right) I_C.
 \end{equation}

 Again first setting $m = n = 0$, we get
 \[
 I_A = \frac{2}{\ell} = \frac{a}{\text{Area}(T)}.
 \]
 Differentiating with respect to $m$ and $n$ we get the two equations
 \begin{equation}
   \label{E:d-m-2}
 0 = I_A + \frac{a_1}{b} I_B - \frac{a + a_1}{c} I_C,
 \end{equation}
 and
 \[
 0 = - \frac{\ell}{b} I_B + \frac{\ell}{c} I_C,
 \]
 so that again
 \[
 I_B = \frac{b}{c} I_C.
 \]
 Now substituting into 
 \eqref{E:d-m-2}, we have
 \[
 \frac{2}{\ell} = \left( - \left(\frac{a_1}{b} \right) \left( \frac{b}{c} \right)  + \frac{a +
   a_1}{c} \right) I_C,
 \]
 or
 \[
 I_C = \frac{2c}{a \ell}  = \frac{c}{\text{Area}(T)}.
 \]
 Back substituting once again, we have also
 \[
 I_B = \frac{b}{\text{Area} (T) }.
 \]
 This completes the proof of the obtuse triangle case, and hence
 proves the theorem.

 \section{Other polygons}
 \label{S:example}
 Theorem \ref{T:main} is false in general for other planar polygons.  One can
 see from the computations above that it is straightforward to come up
 with 3 independent equations relating the Neumann data on the sides.
 The proof suggests that only three equations are possible in
 general.  Of course we do not have a proof of that.  However, even
 for convex polygons the theorem is false in general.

 Consider the square $\Omega = [0, 2 \pi]^2$.  The Dirichlet
 eigenfunctions are given by the Fourier basis.  Let us examine some
 specific choices.  That is, for integers $j,k$ let
 \[
 u_{jk} (x,y) = ( \pi)^{-1} \sin (jx ) \sin(ky).
 \]
 The $u_{jk}$ vanish on $\p \Omega$, and are normalized.  
 We  have
 \[
 -\Delta u_{jk} = (j^2 + k^2 ) u_{jk}
 \]
 as usual.  Rescaling to a semi-classical equation, we take $h = (j^2
 + k^2)^{-1/2}$ to get
 \[
 -h^2 \Delta u_{jk} = u_{jk}.
 \]
 On $x = 0$ and $x = 2 \pi$ repectively, we have
 \[
 h \p_\nu u_{jk}|_{x = 0}  = \pi^{-1} (-j) h \sin(ky)
 \]
 and
 \[
 h \p_\nu u_{jk} |_{x = 2 \pi} = \pi^{-1} j h \sin(ky).
 \]
 Hence the norm of the Neumann data along either $x = 0$ or $x = 2 \pi$ is
 \[
 \int_{0}^{2 \pi} | h \p_\nu u_{jk}|_{x = 0, 2 \pi} |^2 dy = \pi^{-1} h^2 j^2.
 \]
 If $k \gg j$, we can make $h^2 j^2$ as small as we like, so in fact
 we can say only
 \[
 \int_0^{2 \pi} | h \p_\nu u_{jk} |_{x = 0, 2 \pi} |^2 dy \geq c h^2
 \]
 for some $c>0$.  A rough conjecture is that a lower bound of $h^2$
 holds for any polygon in the plane.  The author plans to revisit this
 question in subsequence papers.

\bibliographystyle{alpha}
\bibliography{HC-bib}

\end{document}